\documentclass{article}

\newcommand{\ba}{\begin{array}}
\newcommand{\ea}{\end{array}}
\newcommand{\f}{\varphi}
\usepackage{amssymb}

\begin{document}
\title{\bf Existence theorem on spectral function for singular nonsymmetric first order differential operators}
\author{\large\bf Wuqing Ning \footnote{E-mail: wqning@ustc.edu.cn}}
\date{}

\maketitle
\begin{center}
\large {School of Mathematical Sciences, University of Science and
Technology of China, Hefei 230026, China}
\end{center}

{\bf Abstract}
In this paper we study spectral function for a nonsymmetric
differential operator on the half line. Two cases of the coefficient matrix are considered, and for each case
we prove by Marchenko's method that, to the boundary value problem, there
corresponds a spectral function related to which
a Marchenko-Parseval equality and an expansion formula are established. Our results extend the
classical spectral theory for self-adjoint Sturm-Liouville
operators and Dirac operators.\\



\section{Introduction}
\setcounter{equation}{0}
\quad \  As a very essential mathematical problem, Weyl-Stone eigenfunction expansion \cite{stone,weyl} in which the key role is spectral function for  singular
 self-adjoint second order linear differential operators has been
 studied deeply by many renowned mathematicians such as K. Kodaira \cite{kodaira}, N. Levinson
 \cite{levinson}, B. M. Levitan \cite{lev50},  E. C. Titchmarsh \cite{titchmarsh} and K. Yosida \cite{yosida50}
  and so on.  This theory for the singular operators may be
 derived as a limiting case of the classical Sturm-Liouville
 expansion theorem for the regular operators, where the Parseval equality for the regular
 operators plays a very important role for the proofs. Similar ideas  can also be applied to singular
 self-adjoint first order systems, for example, the Dirac operators \cite{lesch,lev91}.
  For general theory of eigenfunction
 expansion for self-adjoint and regular non-self-adjoint operators in Hilbert space,
 we refer to \cite{berez,kotani,me}. For multidimensional cases, see, e.g., \cite{ike}. Moreover, a two-fold spectral expansion in terms of
principal functions of a Schr\"{o}dinger operator has been derived in \cite{bck}. Recently, K. Kirsten and P. Loya \cite{kl} have
obtained some interesting results on the spectral zeta function
for a Schr\"{o}dinger operator on the half line.

 However, to the author's knowledge, for singular
 nonsymmetric differential operators,  there are few results
on eigenfunction expansion. The limiting approach for self-adjoint case
can not be applied even for very simple case of nonsymmetric
differential operators, since in general the corresponding regular spectrum has irregular behavior on the complex plane.  In order to
 extend expansion theory to general case, V. A. Marchenko \cite{mar63,mar} established an excellent method in
dealing with the singular Sturm-Liouville operator with complex-valued
potential. In this paper, inspired by the idea of V. A. Marchenko, we are going to establish expansion
theorem in two cases for a singular nonsymmetric differential operator, where the key is to prove the existence of the
corresponding spectral function. Our results can be extended to $2n\times 2n$ systems, and for simplicity we here will
only consider the case of $n=1$. For the regular case of this nonsymmetric differential operator, recently  we
have obtained some results on inverse spectral problems with applications to inverse problems
for one-dimensional hyperbolic systems, see \cite{n06}--\cite{ny08}.
It is well known that for many differential operators there are intrinsic relations between their spectral functions and
 the corresponding Weyl functions (often called $m$-functions), and for the recent interesting results on Weyl functions see,
e.g.,  \cite{frei01,gb,gks,ks,sak,yurko92,yurko02}. For the
asymptotic behavior of spectral functions for elliptic operators we
refer to \cite{gording,hormander,miyazaki}.

 In this paper we consider boundary value problems generated
by a nonsymmetric differential
 operator on the half line $0\leq x<\infty$:
 $$
(A_P\f) (x):=B\frac{{\rm d}\f }{{\rm d}
 x}(x)+P(x)\f (x)=\lambda\f(x)
$$
 where  $B=\left (\ba{cc} 0 & 1 \\ 1 & 0 \ea\right )$, both matrix-valued function $P=\left (\ba{ll} p_{11} & p_{12} \\ p_{21} & p_{22} \ea\right )
\in (C^1[0,\infty))^4$ and parameter $\lambda$ are complex-valued. It is directly checked that the adjoint operator of $A_P$
in some suitable Hilbert space is $-B\frac{{\rm d}}{{\rm d}
 x}+\overline{P^T(x)}$ and consequently  $A_P$ is nonsymmetric. Here and henceforth,
$\overline{c}$ denotes the complex conjugate of $c$ and
 $\cdot^T$ denotes the transpose of a vector or matrix under
 consideration. Here we point out that the spectrum problem for $A_P$ with compact matrix-valued function $P$ has been studied  in \cite{tr}.

In order to describe our results properly, we first give some information on distributions and  we refer to \cite{mar} for more details.  Let
${\mathbb{K}}^2(0,\infty)$ denote the set of all square integrable
functions in $(0,\infty)$ with compact support. For $\sigma>0$, we set
${\mathbb{K}}^2_\sigma(0,\infty)=\{f\in{\mathbb{K}}^2(0,\infty):
f(x)=0\ \mbox{for}\ x>\sigma \}$. The entire function $e(\rho)$ is
 called the {\em function of exponential type} if $|e(\rho)|\leq
C\exp(\sigma|{\rm Im}\rho|)$ where the positive constants $C$ and $\sigma$
depend on $e(\rho)$. Moreover, the index
$$\sigma_e=\displaystyle{\overline{\lim_{r\rightarrow\infty}}r^{-1}\ln\left(\max_{|\rho|=r}|e(\rho)|\right)}$$
is called the type of entire function $e(\rho)$. Let linear
topological space $Z$ be the set of all entire exponential type
functions integrable on the real line. The sequence $e_n$
converges to $e$ in $Z$ if
$\lim_{n\rightarrow\infty}\int_{-\infty}^\infty|e_n(\rho)-e(\rho)|{\rm
d}\rho=0$
 and the types $\sigma_n$ of the functions $e_n(\rho)$
are bounded: $\sup \sigma_n<\infty$. The set of all linear
continuous functionals defined on the test space $Z$ will be denoted by
$Z'$ whose components are called {\em distributions} (generalized
functions). The sequence $D_n$ converges to $D$ in $Z'$ if
${\lim_{n\rightarrow\infty}<D_n,e(\rho)>=<D,e(\rho)>}$
for all test functions $e\in Z$.

  In this paper we consider two cases of the coefficient matrix $P$. The first case is special and will be described as follows. Let $P$ be a continuously differentiable
   matrix-valued function  satisfying
 $BP=PB$ and $\mu$ be a complex constant. Here it is easy to see that $P$ is of the form  $\left(\ba{cc} a & b \\
b & a \ea\right)$. Consider the following
boundary value problem
 \begin{equation}\label{eq1}
\left\{\ba{lll}B\displaystyle{\frac{{\rm d}\f}{{\rm d}
 x}(x)+P(x)\f(x)=\lambda\f(x), \ 0<x<\infty,}\\
 \\
\f(0)=\left (\ba{cc} \cosh\mu & \sinh\mu \\ \sinh\mu & \cosh\mu
\ea\right).
 \ea\right.
 \end{equation}
 Let $\f=\f(x,\lambda)$ be the solution to (\ref{eq1}) and
$$\f_{[1]}=\left (\ba{l} \f_{[1]} ^{(1)}
\\
\\ \f_{[1]}^{(2)}\ea \right )\ \mbox{and}\ \f_{[2]}=\left (\ba{l} \f_{[2]}^{(1)}
\\
\\ \f_{[2]}^{(2)}\ea \right )$$
 be the first and the second column vector of the
matrix $\f$, i.e., $\f=(\f_{[1]}\ \ \f_{[2]})$. Similarly we denote the matrix inverse of $\f$ by
$\psi=\f^{-1}=(\psi_{[1]}\ \ \psi_{[2]})$. Now for
$$f=\left (\ba{l} f^{(1)}
\\ f^{(2)}\ea \right )\in \left(L^2(0,\infty)\right)^2,\ \ g
=\left (\ba{l} g^{(1)}
\\ g^{(2)}\ea \right )\in \left(L^2(0,\infty)\right)^2$$
where $\left(L^2(0,\infty)\right)^2$ denotes the product space of $L^2(0,\infty)$, we set
$$\displaystyle{\omega_f^{k}(\rho)=\int_0^\infty f^T(x)\psi_{[k]}(x,i\rho){\rm
d}x},\ \ \displaystyle{\eta_{{g}}^k(\rho)=\int_0^\infty
\f_{[k]}^T(x,i\rho){g(x)}{\rm d}x}\ \ (k=1,2),
$$
where $i=\sqrt{-1}$, $\rho\in{\Bbb R}$. Then we have the first main result of this paper.\\
\\
 {\bf Theorem 1.} {\em It holds for the boundary value problem (\ref{eq1}) that
 \begin{equation}\label{parseval-1}
\int_0^\infty f^T(x)\overline{g(x)}{\rm d}x=\frac{1}{2\pi}\sum_{k=1}^{2}\int_{-\infty}^{\infty}\omega_f^{k}(\rho)
\eta_{\overline{g}}^k(\rho){\rm d}\rho.
\end{equation}
Moreover, for $f\in\left({\mathbb{K}}^2(0,\infty)\right)^2$ with $\omega_f^{k}(\rho),\eta_{{f}}^k(\rho)\in Z\ (k=1,2)$, the following expansion formula holds:
\begin{equation}\label{expansion-1}
\ba{lll} f(x)\displaystyle{=\frac{1}{2\pi}\sum_{k=1}^{2}\int_{-\infty}^{\infty}\omega_f^{k}(\rho)
\f_{[k]}(x,i\rho){\rm d}\rho}
\\
\ \ \ \ \ \
\displaystyle{=\frac{1}{2\pi}\sum_{k=1}^{2}\int_{-\infty}^{\infty}
\eta_f^k(\rho)\psi_{[k]}(x,i\rho){\rm
d}\rho}.\ea
\end{equation}}\\
\\
We often call (\ref{parseval-1}) (or (\ref{parseval-2})) the {\em Marchenko-Parseval equality} which means that a {\em spectral function} exists in corresponding boundary value problem. Historically,
the concept of spectral function came from the classical theory of Weyl.
Theorem 1 implies that $\frac{1}{2\pi}E$ is a spectral function corresponding to problem (\ref{eq1}) with $P$ satisfying $BP=PB$, which is the same as the case of $P=0$.
Here and henceforth $E$ denotes the $2\times 2$ unit matrix.

For general matrix function $P\in{\left(C^1[0,\infty)\right)^4}$ without the constraint $BP=PB$, we also can  show the
 existence of the corresponding  spectral
 function. More precisely, let $Q$ be a $2\times 2$ matrix satisfying $QB+BQ=B$ and
 $Q^2=Q$. It is seen by simple computation that there exists matrix $Q$ satisfying the above conditions, and the
 simplest one is $Q=$diag$(1,0)$. It follows easily from $\det B=1$ that $\det Q=0$. Consider the following
boundary value problems
\begin{equation}\label{eq-phi}
\left\{\ba{lll}B\displaystyle{\frac{{\rm d}\f}{{\rm d}
 x}(x)+P(x)\f(x)=\lambda\f(x), \ 0<x<\infty,}\\
 \\
\f(0)=Q,
 \ea\right.
 \end{equation}
 and
 \begin{equation}\label{eq-phi*}
\left\{\ba{lll}-\displaystyle{\frac{{\rm d}\widetilde{\f}}{{\rm d}
 x}(x)B+\widetilde{\f}(x)P(x)=\lambda\widetilde{\f}(x), \ 0<x<\infty,}\\
 \\
\widetilde{\f}(0)=Q.
 \ea\right.
 \end{equation}
 Denote the solutions to problems (\ref{eq-phi}) and (\ref{eq-phi*}) by $\f(x,\lambda)$ and $\widetilde{\f}(x,\lambda)$, respectively.
 For all $2\times 2$ matrices $f,g\in\left(L^2(0,\infty)\right)^4$, we set
 \begin{equation}\label{def1}
\Phi_f(\rho)=\int_0^\infty f(x) \f(x,i\rho){\rm d}x,\ \
\widetilde{\Phi}_g(\rho)=\int_0^\infty
\widetilde{\f}(x,i\rho)g(x){\rm d}x,
\end{equation}
where $i=\sqrt{-1}$, $\rho\in{\Bbb R}$. Then  we have another main result of this paper.\\
\\
 {\bf Theorem 2.} {\em To the problems (\ref{eq-phi}) and (\ref{eq-phi*}) there corresponds a distribution-valued
 spectral
 function $D=\left(D_{kl}\right)_{1\leq k,l\leq 2}$ such that
 $D=QDQ$, $D_{kl}\in Z'$ and
\begin{equation}\label{parseval-2}
\int_0^\infty f(x)g(x){\rm d}x=\int_{-\infty}^\infty \Phi_f(\rho)
 D(\rho)
 \widetilde{\Phi}_g(\rho){\rm d}\rho.
\end{equation}
Moreover, for $f\in\left({\mathbb{K}}^2(0,\infty)\right)^4$ with $\Phi_f(\rho),\widetilde{\Phi}_f(\rho)\in Z^4$, the following expansion formula holds:
\begin{equation}\label{expansion-2}
f(x)=\int_{-\infty}^\infty \Phi_f(\rho)
 D(\rho)
 \widetilde{\f}(x,i\rho){\rm d}\rho=\int_{-\infty}^\infty \f(x,i\rho)
 D(\rho)
 \widetilde{\Phi}_f(\rho){\rm d}\rho.
 \end{equation}}

Although Theorem 1 and Theorem 2 have shown the existence of spectral function for the singular
nonsymmetric differential operator in two cases,
we here point out that the uniqueness of spectral function for the operator does not hold generally, which is the same as
 that for Sturm-Liouville operators (see, e.g., \cite{lev91}). Moreover, since the spectral function
is distribution-valued, it is not a measure in general, which is different from the case of
self-adjoint Sturm-Liouville operators. Besides, given singular nonsymmetric differential operators with general $P$, it is still
 an open problem to prove the existence of spectral functions under general boundary conditions.
 On the other hand, it is interesting to investigate the corresponding inverse problems,
namely, given spectral functions or Weyl functions, find the differential operators. See
\cite{gelf} for the classical inverse problem to determine the potential of the
Sturm-Liouville operator from its spectral function and \cite{frei12} for determination of singular differential pencils from the Weyl function. Theorem 1 has implied that the uniqueness does not hold generally for the inverse problems,
and we need impose other assumptions for uniqueness. In a forthcoming  paper we will  study the inverse problems for the singular nonsymmetric differential operator.

The paper is composed of four sections. In Section 2 we establish
transformation formulae for our boundary value problems. Sections 3 and 4 are devoted to prove
Theorem 1 and Theorem 2 by transformation formulae, respectively.
\section{Transformation formulae}
\setcounter{equation}{0}
Set
\begin{equation}\label{Omega}
\Omega=\{(x,y)\in{\Bbb R}^2: 0<y<x\}.
\end{equation}
For $P_j=\left(P_{j,kl}\right)_{1\leq k,l\leq
2}\in\left(C^1[0,\infty)\right)^4\ (j=1,2)$, we define
\begin{equation}\label{theta1}
\theta_1(x)=\displaystyle{\frac{1}{2}\int_0^x\left(P_{2,12}+P_{2,21}-P_{1,12}-P_{1,21}\right)(s){\rm
d}s}
\end{equation}
and
\begin{equation}\label{theta2}
\theta_2(x)=\displaystyle{\frac{1}{2}\int_0^x\left(P_{2,11}+P_{2,22}-P_{1,11}-P_{1,22}\right)(s){\rm
d}s}.
\end{equation}
Moreover let us put
\begin{equation}\label{R}
R(P_1,P_2)(x)=\exp\left(-\theta_1(x)\right)\left
(\ba{cc}\cosh\theta_2(x) &
-\sinh\theta_2(x) \\
-\sinh\theta_2(x) & \cosh\theta_2(x) \ea\right ).
\end{equation}
Here we remark that $R(P_1,P_2)(0)=E$,
$R(P_1,P_2)(x)=R^{-1}(P_2,P_1)(x)$ and $R\left
(-\overline{P_1^T},-\overline{P_2^T}\right )(x)=\overline{R\left
(P_2,P_1\right )(x)}$. Let $M_2(\Bbb C)$ be the set of
all $2\times 2$ complex-valued matrices. We first prove the following
lemma.\\
\\
 {\bf Lemma 2.1.} {\em For any $\lambda\in{\Bbb C}$, $Q\in{ M_2(\Bbb C)}$ with $\det Q=0$ and $P_j\in
 \left(C^1[0,\infty)\right)^4$ $(j=1,2)$, let $\f_j=\f_j(x,\lambda)$ satisfy
\begin{equation}\label{eq-Pj}
\left\{ \ba{lll}\displaystyle{B\frac{{\rm d}\f_j(x)}{{\rm d}
 x}+P_j(x)\f_j(x)=\lambda\f_j(x), \ 0< x<\infty,}\\
 \\
\f_j(0)=Q.\ea\right.
 \end{equation}
Then there exists a unique
$K(P_1,P_2;Q)=\left(K_{kl}(P_1,P_2;Q)\right)_{1\leq k,l\leq
2}\in\left(C^1(\overline{\Omega})\right)^4$ independent of $\lambda$
such that for $0\leq x<\infty$ and all $\lambda\in{\Bbb C}$
\begin{equation}\label{transform}
\f_2(x,\lambda)=\displaystyle{R(P_1,P_2)(x)\f_1(x,\lambda)+\int_0^x
{K(P_1,P_2;Q)(x,y)\f_1(y,\lambda)}{\rm d}y}.
\end{equation}
$$\mbox{(transformation formula)}$$
Here $R(P_1,P_2)(x)$ is defined by (\ref{R}).

Moreover, the kernel $K(P_1,P_2;Q)$ is
the unique solution to the following problem of first order system
(\ref{eq-K})$\sim$(\ref{xx-K}):
\begin{equation}\label{eq-K}
\ba{lll} &&\displaystyle{B\frac{\partial K(P_1,P_2;Q)}{\partial
x}(x,y)+\frac{\partial K(P_1,P_2;Q)}{\partial y}(x,y)B} \\
&&\ \
\displaystyle{+P_2(x)K(P_1,P_2;Q)(x,y)-K(P_1,P_2;Q)(x,y)P_1(y)=0, \
(x,y)\in \Omega}.\ea
\end{equation}
\begin{equation}\label{x0-K}
K(P_1,P_2;Q)(x,0)BQ=0\ \ \ \ \quad(0\leq x<\infty).
\end{equation}
\begin{equation}\label{xx-K}
\ba{lll} &&\ \ K(P_1,P_2;Q)(x,x)B-BK(P_1,P_2;Q)(x,x)\\
&&=B\displaystyle{\frac{{\rm d}R(P_1,P_2)}{{\rm
d}x}(x)+P_2(x)R(P_1,P_2)(x)-R(P_1,P_2)(x)P_1(x)}\\
&&\quad\quad\quad\quad\quad\quad\quad\quad\quad(0\leq
x<\infty).\ea
\end{equation}}\\
\\
{\bf Proof.} We prove this lemma by the idea used in \cite{ya88}.
Since $P_j\in\left(C^1[0,\infty)\right)^4$ $(j=1,2)$, it can be
verified directly that, if
$K(P_1,P_2;Q)\in\left(C^1(\overline{\Omega})\right)^4$ is the
unique solution to  problem (\ref{eq-K})$\sim$(\ref{xx-K}), then
(\ref{transform}) holds. Therefore, it is sufficient to prove the
existence and the uniqueness of the solution to  problem
(\ref{eq-K})$\sim$(\ref{xx-K}) for each
$P_1,P_2\in{\left(C^1[0,\infty)\right)^4}$.

For clarity, we reduce the proof to a special case.  By the
condition $\det Q=0$, we may assume that a complex constant $c$
exists such that $q_2=c q_1$ where $q_1,q_2$ are the first column
vector and the second one of $Q$, respectively. Then it is
sufficient to prove the existence and the uniqueness of the solution
to problem (\ref{eq-K})$\sim$(\ref{xx-K}) in the case
$\f_j(0,\lambda)=q_1$, since problem (\ref{eq-Pj}) is linear.
Moreover, since a complex constant
$c^*$ exists such that $q_1=c^*\left(\ba{l}
\cosh\mu\\ \sinh\mu \ea\right)$ where $\mu\in{\Bbb C}$, it can be
reduced to the case $\f_j(0,\lambda)=\left(\ba{l} \cosh\mu\\
\sinh\mu \ea\right)$. In this case, we denote the the solution to
problem (\ref{eq-K})$\sim$(\ref{xx-K}) by $K(P_1,P_2,\mu)(x,y)$, and
(\ref{x0-K}) has the following form:
\begin{equation}\label{x0-K-mu}
\left \{\ba{l}K_{12}(P_1,P_2, \mu)(x,0)=-\tanh \mu\ K_{11}(P_1,P_2,
\mu)(x,0),\\ K_{22}(P_1,P_2, \mu)(x,0)=-\tanh \mu\ K_{21}(P_1,P_2,
\mu)(x,0). \ea \right.
\end{equation}
If we set
\begin{equation}\label{L}
\left \{\ba{l}L_1(x,y)=K_{12}(P_1,P_2,\mu)(x,y)-
K_{21}(P_1,P_2,\mu)(x,y),\\ L_2(x,y)=K_{11}(P_1,P_2,\mu)(x,y)-
K_{22}(P_1,P_2,\mu)(x,y), \\ L_3(x,y)=K_{11}(P_1,P_2,\mu)(x,y)+
K_{22}(P_1,P_2,\mu)(x,y),\\ L_4(x,y)=K_{12}(P_1,P_2,\mu)(x,y)+
K_{21}(P_1,P_2,\mu)(x,y) \ea \right.
\end{equation}
and $L=L(x,y)=\left(L_1(x,y),L_2(x,y),L_3(x,y),L_4(x,y)\right)$, then we can
rewrite (\ref{eq-K})$\sim$(\ref{xx-K}) as follows:
\begin{equation}\label{L12}
\displaystyle{\frac{\partial L_k(x,y)}{\partial x}-\frac{\partial
L_k(x,y)}{\partial y}=f_k(x,y,L)}\ \ \  ((x,y)\in{\Omega},\
k=1,2),
\end{equation}
\begin{equation}\label{L34}
\displaystyle{\frac{\partial L_k(x,y)}{\partial x}+\frac{\partial
L_k(x,y)}{\partial y}=f_k(x,y,L)}\ \ \  ((x,y)\in{\Omega},\
k=3,4),
\end{equation}
\begin{equation}\label{x0-L12}
\ \ \ \ \ \ \ \ \ \ \ \ L_k(x,x)=r_k(x)    \ \ \ \ \ \ \ \ \ (0\leq x<\infty,\ k=1,2),
\end{equation}
\begin{equation}\label{x0-L34}
\left\{\ba{lll}L_3(x,0)=\sinh(2\mu) L_1(x,0)+\cosh(2\mu) L_2(x,0)\\
 L_4(x,0)= -\cosh(2\mu)L_1(x,0)-\sinh(2\mu) L_2(x,0)\ea
\right. \ (0\leq x<\infty),
\end{equation}
where $f_k(x,y,L)={\frac{1}{2}\sum_{m=1}^4}\left(a_{km}(y)+b_{km}(x)\right)L_m(x,y)$
$(1\leq k\leq 4)$, here $a_{km}(y)$, $b_{km}(x)$ $(1\leq k,m\leq
4)$ are linear combinations of two elements of the matrix
functions $P_1(y)$ and $P_2(x)$ respectively, and $r_k\in{C^1[0,\infty)}$ $(k=1,2)$ are dependent only on
$P_1$ and $P_2$.

Integrating (\ref{L12}), (\ref{L34}) with (\ref{x0-L12}) and
(\ref{x0-L34}) along the characteristics $x+y=const.$ and
$x-y=const.$ respectively, we obtain the following integral
equations:
\begin{equation}\label{int-L12}
L_k(x,y)=\displaystyle{\int_y^{\frac{x+y}{2}} f_k(-s+x+y,s,L){\rm
d}s+r_k(\frac{x+y}{2})}\ \  \left((x,y)\in{\overline{\Omega}},
k=1,2\right),
\end{equation}
and
\begin{equation}\label{int-L34}
\ba{lcl} L_k(x,y)&=&\displaystyle{\int_0^y f_k(s+x-y,s,L){\rm d}}s\\
&&+\displaystyle{\int_0^{\frac{x-y}{2}} \left\{\alpha_k
f_1(-s+x-y,s,L)+\beta_k f_2(-s+x-y,s,L)\right\}{\rm
d}s}\\
&&\ \ \ \ \ \ \displaystyle{+\alpha_k r_1(\frac{x-y}{2}) +\beta_k
r_2(\frac{x-y}{2})}
\\
\\
&&\ \ \ \ \ \ \ \ \ \ \ \left((x,y)\in{\overline{\Omega}},
k=3,4\right),\ea
\end{equation}
where $\alpha_3=\sinh(2\mu)$, $\beta_3=\cosh(2\mu)$ and
$\alpha_4=-\cosh(2\mu)$, $\beta_4=-\sinh(2\mu)$.

The unique solution $L\in{\left(C^1(\overline{\Omega})\right)^4}$ to
(\ref{int-L12}) and (\ref{int-L34}) can be obtained by the iteration
method. In fact, setting
$$
L_k^{(0)}(x,y)=0\ \ \ \left((x,y)\in{\overline{\Omega}}, 1\leq k\leq
4\right),
$$
$$
\ba{ccc}L_k^{(n)}(x,y)=\displaystyle{\int_y^{\frac{x+y}{2}}
f_k\left(-s+x+y,s,L^{(n-1)}\right){\rm
d}s+r_k(\frac{x+y}{2})}\\
 \left((x,y)\in{\overline{\Omega}},n\geq
1, k=1,2\right),\ea
$$
and
$$
\ba{lcl} &&L_k^{(n)}(x,y)\\
\\
&=&\displaystyle{\int_0^y f_k\left(s+x-y,s,L^{(n-1)}\right){\rm d}}s\\
\\
&&+\displaystyle{\int_0^{\frac{x-y}{2}} \left\{\alpha_k
f_1\left(-s+x-y,s,L^{(n-1)}\right)+\beta_k
f_2\left(-s+x-y,s,L^{(n-1)}\right)\right\}{\rm
d}s}\\
\\
&&\ \ \ \ \ \ \displaystyle{+\alpha_k r_1(\frac{x-y}{2}) +\beta_k
r_2(\frac{x-y}{2})}
\\
\\
&&\ \ \ \ \ \ \ \ \ \ \ \ \ \ \ \ \ \ \ \ \ \ \ \ \
\left((x,y)\in{\overline{\Omega}}, n\geq 1, k=3,4\right),\ea
$$
we can obtain by induction the estimates for each $n\geq 1$
\begin{equation}\label{iter-est}
\left|L_k^{(n)}(x,y)-L_k^{(n-1)}(x,y)\right|\leq
\omega(x)\frac{\zeta^{n-1}(x)}{(n-1)!}\ \ \
\left((x,y)\in{\overline{\Omega}}, 1\leq k\leq 4\right),
\end{equation}
where
$$\omega(x)=\left(|\sinh(2\mu)|+|\cosh(2\mu)|+1\right)\max_{0\leq s\leq x}\left(|r_1(s)|+|r_2(s)|\right)
$$
and
$$\zeta(x)=\left(|\sinh(2\mu)|+|\cosh(2\mu)|+1\right)\ x \max_{0\leq s\leq x}
\frac{1}{2}\sum_{k,l=1}^{2}\left(|P_{1,kl}(s)|+|P_{2,kl}(s)|\right).
$$
Thus
$L_k(x,y)=\lim_{n\rightarrow\infty}L_k^{(n)}(x,y)$
$(1\leq k\leq 4)$ exist uniformly for
$(x,y)\in{\overline{\Omega}}$ and we see that $L_k(x,y)$ $(1\leq
k\leq 4)$ satisfy (\ref{int-L12}) and (\ref{int-L34}) with the
bound $|L_k(x,y)|\leq \omega(x)\exp(\sigma(x))$.

Moreover, differentiating (\ref{int-L12}) and (\ref{int-L34}) with
respect to $x$ and $y$, we can similarly obtain by induction the following
estimates
\begin{equation}\label{iter-est-px}
\left|\displaystyle{\frac{\partial L_k^{(n)}(x,y)}{\partial
x}-\frac{\partial L_k^{(n-1)}(x,y)}{\partial x}}\right|\leq
\xi(x)\frac{\zeta^{n-1}(x)}{(n-1)!}\ \ \
\left((x,y)\in{\overline{\Omega}}, 1\leq k\leq 4\right),
\end{equation}
\begin{equation}\label{iter-est-py}
\left|\displaystyle{\frac{\partial L_k^{(n)}(x,y)}{\partial
y}-\frac{\partial L_k^{(n-1)}(x,y)}{\partial y}}\right|\leq
\xi(x)\frac{\zeta^{n-1}(x)}{(n-1)!}\ \ \
\left((x,y)\in{\overline{\Omega}}, 1\leq k\leq 4\right),
\end{equation}
where
$$\ba{lll}&\ \ \xi(x)\\
&=\displaystyle{\frac{1}{2}\left(|\sinh(2\mu)|+|\cosh(2\mu)|+1\right)}\\
&\ \ \ \times \displaystyle{\Big\{ \max_{0\leq s\leq
x}\left(|r'_1(s)|+|r'_2(s)|\right)
+\frac{1}{2}\omega(x)\exp(\zeta(x))}\\
&\ \ \ \ \ \ \ \ \ \ \ \ \ \times\displaystyle{\max_{0\leq s\leq
x}\sum_{k,l=1}^{2}\left(|P_{1,kl}(s)|+|P_{2,kl}(s)|+\left(|P'_{1,kl}(s)|+|P'_{2,kl}(s)|\right)x\right)\Big\}.}
\ea
$$
Therefore, it follows from (\ref{iter-est-px}) and
(\ref{iter-est-py}) that
$L\in{\left(C^1(\overline{\Omega})\right)^4}$. The uniqueness of the
solution to (\ref{eq-K})$\sim$(\ref{xx-K}) is shown by
(\ref{iter-est}). $\hfill\square$\\
\\
{\bf Corollary 2.2.} {\em For $j=1,2$, let $\f_j$ be the solution to problem (\ref{eq1}) with $P=P_j\in (C^1[0,\infty))^4$ satisfying
 $P_jB=BP_j$. Then the following
transformation formula holds:
\begin{equation}\label{coro3.2}
\f_2(x,\lambda)=R(P_1,P_2)(x)\f_1(x,\lambda)
\end{equation}
where $R(P_1,P_2)(x)$ is defined by (\ref{R}).}\\

Corollary 2.2 follows from the fact that
$K(P_1,P_2,\mu)\equiv 0$, which can be derived easily by observing that the right hand side of (\ref{xx-K}) is $0$
(in this case the condition $\det Q=0$ is not necessary).
Or one may directly verify
(\ref{coro3.2}). Here we omit the details.\\
\\
{\bf Corollary 2.3.}\\ {\em Let $S$ and $\widetilde{S}$ be the solutions corresponding to
$P=0$ in (\ref{eq-phi}) and (\ref{eq-phi*}), respectively. Then the following transformation formulae hold.\\
 (1) For problem (\ref{eq-phi}) we have
\begin{equation}\label{f2-f1}
S(x,i\rho)=\displaystyle{R(P,0)(x)\f(x,i\rho)+\int_0^x
{K(P,0;Q)(x,y)\f_(y,i\rho)}{\rm d}y}
\end{equation}
where the kernel
$K(P,0;Q)\in\left(C^1(\overline{\Omega})\right)^4$ satisfies the
equation
\begin{equation}\label{K-p1p2}
\displaystyle{B K_x(P,0;Q)(x,y)+K_y(P,0;Q)(x,y)B}-K(P,0;Q)(x,y)P(y)=0,
(x,y)\in \Omega
\end{equation}
as well as the conditions for  $0\leq x<\infty$
\begin{equation}\label{K-p1p2-0}
 K(P,0;Q)(x,0)Q=K(P,0;Q)(x,0),
\end{equation}
and
\begin{equation}\label{K-p1p2-xx}
K(P,0;Q)(x,x)B-BK(P,0;Q)(x,x)=B R'(P,0)(x)-R(P,0)(x)P(x).
\end{equation}
(2) For problem (\ref{eq-phi*}) we have
\begin{equation}\label{f2-f1*}
\widetilde{S}(x,i\rho)=\displaystyle{\widetilde{\f}(x,i\rho)R(0,P)(x)+\int_0^x
{\widetilde{\f}(y,i\rho)}\overline{K^T(-\overline{P^T},0;\overline{Q^T})(x,y)}{\rm
d}y}
\end{equation}
where the kernel
$\overline{K^T(-\overline{P^T},0;\overline{Q^T})(x,y)}$
satisfies
\begin{equation}\label{K-p1p2*}
Q\overline{K^T(-\overline{P^T},0;\overline{Q^T})(x,0)}
=\overline{K^T(-\overline{P^T},0;\overline{Q^T})(x,0)}.
\end{equation}}\\
\\
 {\bf Proof.} (1) is obvious, since $\det Q=0$ and then Lemma 2.1 can be applied. Here (\ref{K-p1p2-0}) follows
from (\ref{x0-K}), $BQ=B-QB$ and $B^2=E$. Now we prove (2).
Note that by (\ref{eq-phi*}) the function $\widetilde{\f}(x,i\rho)$ statifies
$$\left\{\ba{lll}B\displaystyle{\frac{{\rm d}\overline{\widetilde{\f}^T}}{{\rm d}
 x}(x)-\overline{P^T(x)}\overline{\widetilde{\f}^T(x)}=i\rho\overline{\widetilde{\f}^T(x)}, \ 0<x<\infty,}\\
 \\
\overline{\widetilde{\f}^T(0)}=\overline{Q^T}.
 \ea\right.
$$
Then  one will obtain
(\ref{f2-f1*}) by (1) if he notices the following fact:
$R(0,P)(x)=\overline{R(-\overline{P^T},0)(x)}$.$\hfill\square$\\

Since the solutions to the boundary value problems with $P=0$ are entire in $\lambda$, it follows easily from the transformation formulae that\\
\\
{\bf Corollary 2.4.} {\em
For each fixed $x$, all solutions to the boundary value problems under consideration are entire in $\lambda$.}

\section{Proof of Theorem 1}
\setcounter{equation}{0}
We divide the proof of Theorem 1  into four steps as follows.\\
\\
{\em First step.} We first construct a
regular spectral function. Let $S$ denote the solution of (\ref{eq1}) corresponding to $P=0$.
Set
$$\rho=-i\lambda\ \mbox{and}\ \nu=-i\mu.$$
 It is easy to see that
\begin{equation}\label{S}
\ba{lll} S=S(x,\lambda)&=&\left(\ba{cc} \cosh(\lambda x+\mu)& \sinh(\lambda x+\mu)\\
\\ \sinh(\lambda x+\mu) & \cosh(\lambda x+\mu)\ea \right)\\
\\
&=&\left(\ba{cc} \cos(\rho
x+\nu) & i\sin(\rho x+\nu)\\  i\sin(\rho x+\nu) & \cos(\rho
x+\nu)\ea \right)\ea
\end{equation}
and
\begin{equation}\label{S-1}
S^{-1}=S^{-1}(x,\lambda)=\left(\ba{cc} \cos(\rho
x+\nu) & -i\sin(\rho x+\nu)\\  -i\sin(\rho x+\nu) & \cos(\rho
x+\nu)\ea \right).
\end{equation}
 We choose two sufficiently smooth real-valued functions
$\delta_n(x)$ and $\gamma_\sigma(x)$ subject to the following
conditions:
$$\displaystyle{\int_0^\infty\delta_n(x){\rm d}x=1,}
$$
$$\delta_n(x)=0\ {\mbox{for}}\ x=0\ {\mbox{and}}\
x\geq\displaystyle{\frac{1}{n}},\ \ \delta_n(x)>0\ {\mbox{for}}\
0<x< \displaystyle{\frac{1}{n}},
$$
\begin{equation}\label{gamma-sigma}
\gamma_\sigma(x)=1\ {\mbox{for}}\ \
0\leq x\leq \sigma,\ \ \gamma_\sigma(x)=0\ {\mbox{for}}\ x>
\sigma+1,
\end{equation}
and it is obvious that $\delta_n(x)$ tends to the Dirac delta function $\delta(x)$
as $n\rightarrow\infty$. We set
\begin{equation}\label{Dn}
\ba{lll}&D_n^\sigma(\rho)=\left(D_{n,jm}^\sigma(\rho)\right)_{1\leq j,m\leq 2}\\
\\ &\ \ \ \ \ \ \ \ \ =\displaystyle{\frac{1}{2\pi}\int_{0}^{\infty}
\left(\ba{cc} \cos(\rho x+\nu) & -i\sin(\rho x+\nu) \\
-i\sin(\rho x+\nu) & \cos(\rho x+\nu)\ea\right)}\\
\\
&\ \ \ \ \ \ \ \ \ \ \ \ \times R(P,0)(x)\delta_n(x) E \gamma_\sigma(x)
\left(\ba{cc} \cos\nu & i\sin\nu \\
i\sin\nu & \cos\nu\ea\right){\rm d}x.\ea
\end{equation}
Since the Fourier transform is a one-to-one mapping on the space of bounded continuous
Lebesgue-integrable functions
and $R(P,0)(x)\delta_n(x) E \gamma_\sigma(x)$ is
a continuously differentiable matrix function with
compact support, it is not hard to see that the matrix
function $D_n^\sigma(\rho)$ is bounded and Lebesgue-integrable on the real line
$-\infty<\rho<\infty$. Hence the integral
$$\displaystyle{\int_{-\infty}^{\infty}
S(x,i\rho)
D_n^\sigma(\rho) \left(\ba{cc} \cos\nu & -i\sin\nu \\
-i\sin\nu & \cos\nu\ea\right){\rm d}\rho}
$$
converges absolutely. By Corollary 2.2 we have
$\f(x,i\rho)$$=R(0,P)(x)S(x,i\rho)$$=R^{-1}(P,0)(x)S(x,i\rho)$, which implies by
the Fourier inverse transform that
\begin{equation}\label{inv-f}
\displaystyle{\int_{-\infty}^{\infty}
\f(x,i\rho)
D_n^\sigma(\rho) \left(\ba{cc} \cos\nu & -i\sin\nu \\
-i\sin\nu & \cos\nu\ea\right){\rm d}\rho=\delta_n(x) E\ \ (0\leq
x\leq \sigma).}
\end{equation}
Here and henceforth we repeatedly
make use of the fact that two matrices $P_1$ and $P_2$ in the form of $\left(\ba{cc} a & b \\
b & a \ea\right)$ are interchangeable: $P_1P_2=P_2P_1$.\\
\\
{\em Second step.} Next we will investigate the asymptotic behavior of the following matrix function
as $n\rightarrow\infty$
\begin{equation}\label{eq-Un-sigma}
\ba{ll} U_n^\sigma(x,y)=\left(U_{n,kl}^\sigma(x,y)\right)_{1\leq k,l\leq 2}\\
\\
\quad\quad\quad\ \  :=\displaystyle{\int_{-\infty}^{\infty}
\f(x,i\rho) D_n^\sigma(\rho) \f^{-1}(y,i\rho){\rm d}\rho\ \ \ (0\leq x,y\leq \sigma).}
\ea
\end{equation}
It is easy to find that
$$U_n^\sigma(x,0)=\displaystyle{\int_{-\infty}^{\infty}
\f(x,i\rho) D_n^\sigma(\rho) \left(\ba{cc} \cos\nu & -i\sin\nu \\
-i\sin\nu & \cos\nu\ea\right) {\rm d}\rho=\delta_n(x) E\ \ (0\leq
x\leq \sigma)}
$$
and
$$U_n^\sigma(0,y)=\displaystyle{\int_{-\infty}^{\infty}
 \left(\ba{cc} \cos\nu & i\sin\nu \\
i\sin\nu & \cos\nu\ea\right)D_n^\sigma(\rho)\f^{-1}(y,i\rho)  {\rm
d}\rho\ \ (0\leq y\leq \sigma).}
$$
Now we show that $U_n^\sigma(0,y)=0$ for all $y\geq 0$. Indeed,
first one can see from (\ref{Dn}) that
$$\ba{lll}D_n^\sigma(\rho)&=&\displaystyle{\frac{1}{2\pi}\int_{0}^{\infty}
\left(\ba{cc} \cos(\rho x) & -i\sin(\rho x) \\
-i\sin(\rho x) & \cos(\rho x)\ea\right)}\\
\\
&&\ \ \ \ \ \ \times \left\{
R_{11}(P,0)(x)E+R_{12}(P,0)(x)B\right\}\delta_n(x) \gamma_\sigma(x)
{\rm d}x.\ea
$$
Moreover, for any continuous
scalar function $u(x)$ with compact support and $u(0)=0$, it follows easily from the theory of
the Fourier cosine and sine transforms that
\begin{equation}\label{c-s}
\ba{ll} \ \ \ \displaystyle{\int_{-\infty}^{\infty}\int_0^{\infty}
 \left(\ba{cc} \cos(\rho x) & -i\sin(\rho x) \\
-i\sin(\rho x) & \cos(\rho x)\ea\right)u(x)\left(\ba{cc} \cos(\rho y) & -i\sin(\rho y) \\
-i\sin(\rho y) & \cos(\rho y)\ea\right) {\rm d}x {\rm d}\rho}\\
\\
=0.\ea
\end{equation}
 Consequently, it follows from (\ref{c-s}) and $\f^{-1}(\cdot,i\rho)=S^{-1}(\cdot,i\rho)R(P,0)(\cdot)$ that
 $U_n^\sigma(0,y)R(0,P)(y)=0$ and hence $U_n^\sigma(0,y)=0$, since $R(0,P)(y)$ is invertible.

 On the other hand,  by (\ref{eq-Un-sigma}) it is easy to see that, for fixed $n$ and $\sigma$, $U_n^\sigma(\sigma,\cdot)$ is
  a bounded differentiable function on $[0,\sigma]$ and denoted by $\Xi_n(\cdot)$ for simplicity. Therefore, since by (\ref{eq1}) we easily show that
$${\frac{{\rm d}\f^{-1}(x)}{{\rm
d}x}B-\f^{-1}(x)P(x)=-i\rho\f^{-1}(x)},
$$
the above argument implies that the functions
\begin{equation}\label{UnN}
U_{nN}^\sigma(x,y):=\displaystyle{ \int_{-N}^{N} \f(x,i\rho)
D_n^\sigma(\rho) \f^{-1}(y,i\rho){\rm d}\rho}
\end{equation}
are continuously differentiable and satisfy the equation
\begin{equation}\label{eq-Un-N}
 \displaystyle{B\frac{\partial U}{\partial
x}(x,y)+\frac{\partial U}{\partial y}(x,y)B}
\displaystyle{+P(x)U(x,y)-U(x,y)P(y)=0 \ \mbox{in}\ \Pi_{\sigma}}
\end{equation}
 as well as the following conditions
\begin{equation}\label{bc-Un-N}
U(x,0)=\delta_{nN}(x) E\ \ \ (0\leq x\leq\sigma),
\end{equation}
\begin{equation}\label{bc-Un-N-2}
U(0,y)=\Gamma_{nN}(y),\ \ U(\sigma,y)=\Xi_{nN}(y) \ \ \ (0\leq y\leq\sigma),
\end{equation}
where $\Pi_{\sigma}=\{(x,y)\in{\Bbb R}^2: 0<x,y<\sigma\}$, the functions $\delta_{nN},\Gamma_{nN}$ and $\Xi_{nN}$
satisfy the compatibility conditions and
${\lim_{N\rightarrow\infty}\delta_{nN}(x)}$$=\delta_n(x)$,
${\lim_{N\rightarrow\infty}\Gamma_{nN}(y)}=0$ and $\lim_{N\rightarrow\infty}\Xi_{nN}(y)=\Xi_n(y)$.
We should note
that problem (\ref{eq-Un-N}), (\ref{bc-Un-N}) and (\ref{bc-Un-N-2})
can be rewritten as a symmetric hyperbolic
system:\\
\begin{equation}\label{hyper-V}
\left\{\ba{lll}\displaystyle{\frac{\partial V}{\partial y}(x,y)+\left(\ba{cc} 0 & E \\
E & 0\ea\right)\frac{\partial V}{\partial x}(x,y)+C(x,y)V(x,y)=0}\ \mbox{in}\ \Pi_{\sigma},\\
\\
V(x,0)=\delta_{nN}(x) \overrightarrow{H}\ \ \ \  (0\leq x\leq\sigma),\\
\\
V(0,y)=\overrightarrow{\Gamma}_{nN}(y),\ \ V(\sigma,y)=\overrightarrow{\Xi}_{nN}(y) \ \ \ (0\leq y\leq\sigma),\ea\right.
\end{equation}
where
$$V(x,y)=\left(\ba{cc} U_{11}(x,y)\\U_{12}(x,y)\\U_{21}(x,y)\\U_{22}(x,y)
\ea\right),\ \ \overrightarrow{H}=\left(\ba{cc} 1\\0\\0\\1 \ea\right),$$
$$\overrightarrow{\Gamma}_{nN}(y)=\left(\ba{cc}
\Gamma_{nN,11}(y)\\ \Gamma_{nN,12}(y)\\ \Gamma_{nN,21}(y)\\ \Gamma_{nN,22}(y)
\ea\right),\ \ \overrightarrow{\Xi}_{nN}(y)=\left(\ba{cc}
\Xi_{nN,11}(y)\\ \Xi_{nN,12}(y)\\ \Xi_{nN,21}(y)\\ \Xi_{nN,22}(y)
\ea\right)
$$
and $C(x,y)$ is the following $4\times 4$ matrix-valued function
$$\left(\ba{cccc} -P_{12}(y) & P_{11}(x)-P_{22}(y) & 0 & P_{12}(x) \\
\\
P_{11}(x)-P_{11}(y) & -P_{21}(y) & P_{12}(x) & 0 \\
\\
0 & P_{21}(x) & -P_{12}(y) & P_{22}(x)-P_{22}(y) \\
\\
P_{21}(x) & 0 & P_{22}(x)-P_{11}(y) & -P_{21}(y)\ea\right).
$$
Since $BU_{nN}^\sigma(x,y)=U_{nN}^\sigma(x,y)B$, a direct calculation shows that the symmetric hyperbolic
system (\ref{hyper-V}) is actually equivalent to the following normal hyperbolic system
\begin{equation}\label{eq-v}
\left\{\ba{lll}\displaystyle{\frac{\partial v}{\partial y}(x,y)=\left(\ba{cc} 1 & 0 \\
0 & -1 \ea\right)\frac{\partial v}{\partial x}(x,y)+c(x,y)v(x,y)}\ \mbox{in}\ \Pi_{\sigma},\\
\\
v(x,0)=\delta_{nN}(x) \overrightarrow{h}\ \  (0\leq x\leq\sigma),\\
\\
v_2(0,y)=2v_1(0,y)-\Gamma_{nN,11}(y)+3\Gamma_{nN,12}(y),\\
\\
v_2(\sigma,y)=2v_1(\sigma,y)-\Xi_{nN,11}(y)+3\Xi_{nN,12}(y) \ \ \ (0\leq y\leq\sigma),\ea\right.
\end{equation}
where
$$v(x,y)=\left(\ba{cc} v_1(x,y)\\ v_2(x,y) \ea\right)=\left(\ba{cc} U_{11}(x,y)-U_{12}(x,y)\\ U_{11}(x,y)+U_{12}(x,y) \ea\right),\ \overrightarrow{h}= \left(\ba{cc} 1\\ 1 \ea\right),
$$
and
$$\ba{ll}
\ \ c(x,y)\\
\\
=\left(\ba{cc} (P_{11}-P_{12})(x)+(P_{12}-P_{11})(y)  &  (P_{12}-P_{11})(x)+(P_{22}-P_{21})(y)\\
\\
(P_{11}+P_{12})(y)-(P_{11}+P_{12})(x)  & (P_{22}+P_{21})(y)-(P_{11}+P_{12})(x)
\ea\right).
\ea
$$
If we take the variable $y$ as time, then it is not hard to verify that the classical
Uniform Kreiss Condition holds, and hence from the well-known
results of well-posedness on linear hyperbolic systems (cf. \cite{higdon} and references therein) we see
that (\ref{eq-v}) has a unique solution, that is, there exists a
unique solution $U_{nN}^\sigma(x,y)$ to problem (\ref{eq-Un-N}),
(\ref{bc-Un-N}) and (\ref{bc-Un-N-2}) such that
$U_{nN}^\sigma(x,y)\rightarrow U_n^\sigma(x,y)$ as
$N\rightarrow\infty$.

 On the other
hand, if we set
$W_{nN}^\sigma(x,y)=U_{nN}^\sigma(x,y)-\delta_{nN}(x-y) E$ for $0\leq
x,y\leq\sigma$ where $\delta_{nN}(x-y)=0$ for $0\leq
x<y\leq\sigma$, then $W_{nN}^\sigma(x,y)$ satisfies the following equation
\begin{equation}\label{hyper-W}
\ba{ll}&\ \ \displaystyle{ B\frac{\partial W}{\partial
x}(x,y)+\frac{\partial W}{\partial y}(x,y)B
+P(x)W(x,y)-W(x,y)P(y)}\\
\\
&=\delta_{nN}(x-y)\left(P(y)-P(x)\right)
 \ea
\end{equation}
and $W(x,0)=0$. It follows easily from the compatibility conditions that $W_{nN}^\sigma(0,y)\rightarrow 0$ $(N,n\rightarrow \infty)$. Next we will show
\begin{equation}\label{w-sigma}
W_{nN}^\sigma(\sigma,y)\rightarrow 0 \ \ (N,n\rightarrow \infty).
\end{equation}
In fact, it follows from (\ref{S}), (\ref{S-1}), (\ref{inv-f}) and the transformation formulae
 $\f(\cdot,i\rho)$$=R(0,P)(\cdot)S(\cdot,i\rho)$$=R^{-1}(P,0)(\cdot)S(\cdot,i\rho)$, $\f^{-1}(\cdot,i\rho)=S^{-1}(\cdot,i\rho)R(P,0)(\cdot)$  that
$$\ba{ll}
\quad \Xi_n(y)\\
\\
\displaystyle{=R(0,P)(\sigma)R(P,0)(y)R(P,0)(\sigma-y)}\\
\\
\displaystyle{\quad\ \times\int_{-\infty}^{\infty}
\f(\sigma-y,i\rho)
D_n^\sigma(\rho) \left(\ba{cc} \cos\nu & -i\sin\nu \\
-i\sin\nu & \cos\nu\ea\right){\rm d}\rho}\\
\\
=R(0,P)(\sigma)R(P,0)(y)R(P,0)(\sigma-y)\delta_n(\sigma-y) E\ \ (0\leq
x\leq \sigma).
\ea
$$
Thus, we have
\begin{equation}\label{key}
\ba{ll}
\quad\Xi_n(y)-\delta_n(\sigma-y) E\\
\\
=\delta_n(\sigma-y)R(0,P)(\sigma)[R(P,0)(y)R(P,0)(\sigma-y)-R(P,0)(\sigma)]
\ea
\end{equation}
whence (\ref{w-sigma}) follows easily. Consequently, by the well-posedness of symmetric hyperbolic linear
differential equations, we have
$$W_{nN}^\sigma(x,y)\rightarrow 0\ \mbox{as}\ N,n\rightarrow \infty
$$
since $\delta_{nN}(x-y)\rightarrow \delta(x-y)\ \mbox{as}\
N,n\rightarrow \infty$ and hence the right hand side of
(\ref{hyper-W}) tends to $0$.
Therefore, for $0\leq x,y\leq\sigma$
\begin{equation}\label{lim-delta}
U_n^\sigma(x,y)\rightarrow\delta(x-y) E\ \ (n\rightarrow\infty).
\end{equation}
{\em Remark.} There is another and simpler way to prove (\ref{lim-delta}) in which it is not needed to consider (\ref{eq-Un-N}).
The key idea is based on considering (\ref{eq-Un-sigma}), (\ref{c-s}) and (\ref{key}) with replacing $\sigma$ by $x$. We leave
the details to the reader.\\
\\
{\em Third step.} We prove the Marchenko-Parseval equality (\ref{parseval-1}). Assuming that
$f,g\in{\left({\mathbb{K}}^2_\sigma(0,\infty)\right)^2}$
 have compact support, we have by changing the
order of integration that
$$\ba{lll}
&\ \
\displaystyle{\int_0^\infty\int_0^\infty\sum_{k,l=1}^{2}U_{n,kl}^\sigma(x,y)\overline{g^{(k)}(x)}f^{(l)}(y){\rm
d}x{\rm d}y}\\
\\
&=\displaystyle{\int_0^\infty\int_0^\infty\sum_{k,l=1}^{2}\left(\int_{-\infty}^{\infty}\sum_{j,m=1}^{2}
D_{n,jm}^\sigma(\rho) \f_{[j]}^{(k)}(x,i\rho)
\psi_{[m]}^{(l)}(y,i\rho){\rm d}\rho\right)}\\
\\
&\ \ \ \ \ \ \ \ \ \ \ \ \ \ \ \ \ \ \ \ \ \ \ \ \times
\overline{g^{(k)}(x)}f^{(l)}(y){\rm d}x{\rm d}y\\
\\
&=\displaystyle{\sum_{j,m=1}^{2}\int_{-\infty}^{\infty}{\rm
d}\rho\
D_{n,jm}^\sigma(\rho)}\\
\\
&\ \ \ \ \ \ \ \ \ \times
\displaystyle{\left(\int_0^\infty\int_0^\infty\sum_{k,l=1}^{2}
f^{(l)}(y)
\psi_{[m]}^{(l)}(y,i\rho)\f_{[j]}^{(k)}(x,i\rho)\overline{g^{(k)}(x)}{\rm
d}x{\rm d}y\right) }\\
\\
&=\displaystyle{\sum_{j,m=1}^{2}\int_{-\infty}^{\infty}D_{n,jm}^\sigma(\rho)\omega_f^{m}(\rho)
\eta_{\overline{g}}^j(\rho){\rm d}\rho}.
 \ea
$$
Therefore, in view of (\ref{lim-delta}), we
obtain by letting $n\rightarrow\infty$ that for any $f,g\in{\left({\mathbb{K}}^2_\sigma(0,\infty)\right)^2}$
\begin{equation}\label{Par-1}
\displaystyle{\int_0^\infty f^T(x)\overline{g(x)}{\rm d}x=\lim_{n\rightarrow\infty}\sum_{j,m=1}^{2}\int_{-\infty}^{\infty}D_{n,jm}^\sigma(\rho)\omega_f^{m}(\rho)
\eta_{\overline{g}}^j(\rho){\rm d}\rho}.
\end{equation}
By the definition of $D_n^\sigma(\rho)$ (see (\ref{Dn})), we easily see  that
\begin{equation}\label{lim-1}
\lim_{n\rightarrow\infty}D_n^\sigma(\rho)=\frac{1}{2\pi}R(P,0)(0)=\frac{1}{2\pi}E.
\end{equation}
On the other hand, since the Fourier transform is a continuous mapping of $L^2({\Bbb R})$ into $L^2({\Bbb R})$, it follows easily from Corollary 2.2
and the zero extensions of $f$ and $g$ on ${\Bbb R}$ that both $\omega_f^{m}(\rho)$ and $\eta_{\overline{g}}^j(\rho)$ belong to $L^2({\Bbb R})$. Therefore,
Combining (\ref{Par-1}) and (\ref{lim-1}) and letting $\sigma\rightarrow \infty$, we can assert (\ref{parseval-1}) by the boundedness  of $D_n^\sigma(\cdot)$, the dominated convergence theorem and the fact that $\left({\mathbb{K}}^2(0,\infty)\right)^2$ is dense in $\left(L^2(0,\infty)\right)^2$.\\
\\
{\em Forth step.} We prove the expansion (\ref{expansion-1}).
First we assume that $f\in{\left(C_0[0,\infty)\right)^2}$, where $\left(C_0[0,\infty)\right)^2$ denotes the product space of the set of all continuous functions with compact support. For any fixed real number $x\geq 0$ and $\delta>0$, set
\begin{equation}\label{varsigma}
\varsigma(t)=\left\{\ba{ll}\displaystyle{\frac{1}{\delta}}\ \
\mbox{for}\
 t\in (x,x+\delta),\\
 \\
0 \ \ \mbox{for other case} .\ea \right.
\end{equation}
In (\ref{parseval-1})  first letting
$g^{(1)}(t)=\varsigma(t),g^{(2)}(t)=0$ and then letting
$g^{(1)}(t)=0,g^{(2)}(t)=\varsigma(t)$, we have
$$\displaystyle{\frac{1}{\delta}\int_x^{x+\delta}f(t){\rm d}t=
\frac{1}{2\pi}\sum_{k=1}^{2}\int_{-\infty}^{\infty}\omega_f^{k}(\rho)
\frac{1}{\delta} \int_x^{x+\delta} \f_{[k]}(t,i\rho){\rm d}t{\rm
d}\rho}.
$$
Since
$$\displaystyle{\lim_{\delta\rightarrow 0}\frac{1}{\delta}\int_x^{x+\delta}f(t){\rm d}t=f(x)}
$$
and  in $Z$
$$
\displaystyle{\lim_{\delta\rightarrow
0}\omega_f^{k}(\rho)\frac{1}{\delta} \int_x^{x+\delta}
\f_{[k]}(t,i\rho){\rm d}t=\omega_f^{k}(\rho)\f_{[k]}(x,i\rho)},
$$
we prove the first part of (\ref{expansion-1}) by the dominated convergence theorem if
$f\in{\left(C_0[0,\infty)\right)^2}$. For the case of $f\in\left({\mathbb{K}}^2(0,\infty)\right)^2$ we can
approximate $f$ by the functions in
$\left(C_0[0,\infty)\right)^2$. The second part of
(\ref{expansion-1}) can be proved similarly.

\section{Proof of Theorem 2}
First let us prove Theorem 2 for a special case. Recall that $S$ and $\widetilde{S}$ are the solutions corresponding to
$P=0$ in (\ref{eq-phi}) and (\ref{eq-phi*}), respectively.\\
\\
{\bf Lemma 4.1.} {\em For $f,g\in{\left(L^2(0,\infty)\right)^4}$,
it holds that
$$\int_0^\infty f(x)g(x){\rm
d}x=\frac{1}{\pi}\int_{-\infty}^\infty
\Theta_f(\rho)\widetilde{\Theta}_g(\rho){\rm
d}\rho=\frac{1}{\pi}\int_{-\infty}^\infty
\Theta_f(\rho)Q\widetilde{\Theta}_g(\rho){\rm d}\rho$$ and for $x>0$
$$f(x)=\frac{1}{\pi}\int_{-\infty}^\infty
\Theta_f(\rho)\widetilde{S}(x,i\rho){\rm d}\rho=
\frac{1}{\pi}\int_{-\infty}^\infty
S(x,i\rho)\widetilde{\Theta}_f(\rho){\rm d}\rho,$$
where $\Theta_f(\rho)$ and $\widetilde{\Theta}_g(\rho)$ are defined by
\begin{equation}\label{Theta}
\Theta_f(\rho)=\int_0^\infty f(x)S(x,i\rho){\rm d}x,\ \  \widetilde{\Theta}_g(\rho)=\int_0^\infty \widetilde{S}(x,i\rho) g(x){\rm
d}x.
\end{equation}
}\\
 {\bf Proof.} Since it is easy to find that
\begin{equation}\label{SS}
S(x,i\rho)=Q\cosh(i\rho x)+BQ\sinh(i\rho x),\ \widetilde{S}(x,i\rho)=Q\cosh(i\rho x)-QB\sinh(i\rho
 x),
\end{equation}
we have
$$\Theta_f(\rho)=\int_0^\infty f(x)S(x,i\rho){\rm
d}x=\frac{1}{2}\widehat{f}(\rho)(Q-BQ)+\frac{1}{2}\widehat{f}(-\rho)(Q+BQ),
$$
$$\widetilde{\Theta}_g(\rho)=\int_0^\infty \widetilde{S}(x,i\rho) g(x){\rm
d}x=\frac{1}{2}(Q+QB)\widehat{g}(\rho)+\frac{1}{2}(Q-QB)\widehat{g}(-\rho),
$$
where $\widehat{f}(\rho)=\int_0^\infty f(x)\exp(-i\rho x){\rm d}x$ denotes the Fourier transform of $f(x)$.
Therefore, by the well-known Parseval equality
$$\displaystyle{\int_0^\infty f(x)g(x){\rm
d}x=\frac{1}{2\pi} \int_{-\infty}^\infty
\widehat{f}(\rho)\widehat{g}(-\rho){\rm d}\rho}$$
 and the identity for $u,v\in{L^2(0,\infty)}$
 $$\int_{-\infty}^\infty
\widehat{u}(\rho)\widehat{v}(\rho){\rm d}\rho=0,
 $$
 we easily obtain  (note that $Q^2=Q$)
$$\frac{1}{\pi}\int_{-\infty}^\infty
\Theta_f(\rho)\widetilde{\Theta}_g(\rho){\rm
d}\rho=\frac{1}{\pi}\int_{-\infty}^\infty
\Theta_f(\rho)Q\widetilde{\Theta}_g(\rho){\rm d}\rho=\int_0^\infty
f(x)g(x){\rm d}x.$$
On the other hand, since for all $u\in{L^2(0,\infty)}$ and $x>0$ it
holds that
$$\int_{-\infty}^\infty\widehat{u}(\rho)\exp(-i\rho x){\rm d}\rho
=\int_{-\infty}^\infty\widehat{u}(-\rho)\exp(i\rho x){\rm d}\rho=0,
$$
we have
$$\ba{lll}&&\displaystyle{\frac{1}{\pi}\int_{-\infty}^\infty
\Theta_f(\rho)\widetilde{S}(x,i\rho){\rm
d}\rho}\\
\\
&=&\displaystyle{\frac{1}{2}f(x)(Q-BQ)(Q-QB)+\frac{1}{2}f(x)(Q+BQ)(Q+QB)=f(x).}
\ea
$$
 Similarly, we can show that
$\displaystyle{\frac{1}{\pi}\int_{-\infty}^\infty
S(x,i\rho)\widetilde{S}_f(\rho){\rm d}\rho=f(x)}$. $\hfill\square$\\
\\
If we put
\begin{equation}\label{fF}
\displaystyle{f(x)=F(x)R(P,0)(x)+\int_x^\infty
F(t)K(P,0;Q)(t,x){\rm d}t}
\end{equation}
 and
\begin{equation}\label{gG}
\displaystyle{g(x)=R(0,P)(x)G(x)+\int_x^\infty
\overline{K^T\left(-\overline{P^T},0;\overline{Q^T}\right)(t,x)}G(t){\rm
d}t},
\end{equation}
where $F$ and $G$ can be obtained by solving the above Volterra equations of the second kind,
then it follows from changing the order of integration and the transformation formulae
(\ref{f2-f1}) and (\ref{f2-f1*}) that
\begin{equation}\label{fFgG}\Phi_f(\rho)=\Theta_F(\rho),\
\widetilde{\Phi}_g(\rho)=\widetilde{\Theta}_G(\rho).
\end{equation}
Furthermore, we have\\
\\
{\bf Lemma 4.2.} {\em For $f,g\in{\left(L^2(0,\infty)\right)^4}$, it holds that
\begin{equation}\label{fg}\int_0^\infty f(x)g(x){\rm d}x=\int_0^\infty F(x)G(x){\rm d}x+\int_0^\infty\int_0^\infty
F(y) {\mathfrak{F}}(x,y)G(x){\rm d}x{\rm d}y,
\end{equation}
where ${\mathfrak{F}}(x,y)$ is defined as follows:\\
\begin{equation}\label{Fxy}
\ba{ll}\ \ {\mathfrak{F}}(x,y)
=\left\{\ba{lll}
\displaystyle{R(P,0)(y)\overline{K^T\left(-\overline{P^T},0;\overline{Q^T}\right)(x,y)}}\\
\ \ \ \displaystyle{+\int_0^y
K(P,0;Q)(y,t)\overline{K^T\left(-\overline{P^T},0;\overline{Q^T}\right)(x,t)}{\rm
d}t,\ 0\leq y\leq
x,}\\
\\ \displaystyle{ K(P,0;Q)(y,x)R(0,P)(x)}\\
\ \ \ \displaystyle{+\int_0^x
K(P,0;Q)(y,t)\overline{K^T\left(-\overline{P^T},0;\overline{Q^T}\right)(x,t)}
{\rm d}t,\ 0\leq x\leq y.} \ea\right. \ea
\end{equation}}
{\bf Proof.} On one hand, since
$R(P,0)(\cdot)=R^{-1}(0,P)(\cdot)$, we have by changing of the order
of integration
$$\ba{lll}&\ \ \displaystyle{\int_0^\infty f(x)g(x){\rm d}x}\\
\\
&= \displaystyle{\int_0^\infty \left\{ F(x)R(P,0)(x)+\int_x^\infty
F(t)K(P,0;Q)(t,x){\rm d}t\right\}}\\
\\& \ \ \ \ \ \ \ \ \ \  \displaystyle{\times \left\{
R(0,P)(x)G(x)+\int_x^\infty
\overline{K^T\left(-\overline{P^T},0;\overline{Q^T}\right)(t,x)}G(t){\rm
d}t\right\}{\rm d}x}\\
\\
&=\displaystyle{\int_0^\infty F(x)G(x){\rm
d}x+\int_0^\infty\int_x^\infty
F(x)R(P,0)(x)\overline{K^T\left(-\overline{P^T},0;\overline{Q^T}\right)(t,x)}G(t)
{\rm d}t {\rm d}x }
\\
\\
&\ \ \ \ \ +\displaystyle{\int_0^\infty\int_0^t
F(t)K(P,0;Q)(t,x)R(0,P)(x)G(x){\rm d}x {\rm d}t}\\
\\
&\ \ \ \ +\displaystyle{\int_0^\infty\int_x^\infty\int_x^\infty
F(t)K(P,0;Q)(t,x)\overline{K^T\left(-\overline{P^T},0;\overline{Q^T}\right)(s,x)}G(s)
{\rm d}t{\rm d}s{\rm d}x}.\ea$$
 On the other hand, by (\ref{Fxy}),
$$\ba{lll}&\ \  \displaystyle{\int_0^\infty\int_0^\infty
F(y) {\mathfrak{F}}(x,y)G(x){\rm d}x{\rm d}y}\\
\\ &=\displaystyle{\int_0^\infty\int_0^y
F(y) \Big\{ K(P,0;Q)(y,x)R(0,P)(x)}\\
\\ & \ \ \ \ \ \ \ \ \ \ \ \ +\displaystyle{\int_0^x
K(P,0;Q)(y,t)\overline{K^T\left(-\overline{P^T},0;\overline{Q^T}\right)(x,t)}
{\rm d}t \Big\} G(x){\rm d}x{\rm d}y}
\\
\\ &\ \ +\displaystyle{\int_0^\infty\int_y^\infty
F(y) \Big\{
R(P,0)(y)\overline{K^T\left(-\overline{P^T},0;\overline{Q^T}\right)(x,y)}}\\
\\
&\ \ \ \ \ \ \ \ \ \ \ \ +\displaystyle{\int_0^y
K(P,0;Q)(y,t)\overline{K^T\left(-\overline{P^T},0;\overline{Q^T}\right)(x,t)}{\rm
d}t  \Big\}G(x){\rm d}x{\rm d}y}. \ea
$$
Therefore, to prove (\ref{fg}), it is equivalent to show
$$\ba{lll}&\ \ \ \displaystyle{\int_0^\infty\int_x^\infty\int_x^\infty
F(t)K(P,0;Q)(t,x)\overline{K^T\left(-\overline{P^T},0;\overline{Q^T}\right)(s,x)}G(s)
{\rm d}t{\rm d}s{\rm d}x}\\
\\&=\displaystyle{\int_0^\infty\int_0^y\int_0^x F(y)
K(P,0;Q)(y,t)\overline{K^T\left(-\overline{P^T},0;\overline{Q^T}\right)(x,t)}
 G(x){\rm d}t {\rm d}x{\rm d}y}\\
\\ &\ \ \ +\displaystyle{\int_0^\infty\int_y^\infty
 \int_0^y F(y)
K(P,0;Q)(y,t)\overline{K^T\left(-\overline{P^T},0;\overline{Q^T}\right)(x,t)}
G(x){\rm d}t {\rm d}x{\rm d}y },\ea
$$
which can be easily proved by changing of the order of integration.
$\hfill\square$\\
\\
{\em Remark.} It follows easily from (\ref{R}) and Lemma 2.1 that ${\mathfrak{F}}(\cdot,\cdot)\in\left(C^1(\overline{\Omega})\right)^4$ and ${\mathfrak{F}}(\cdot,\cdot)\in\left(C^1(\overline{{\Bbb {R}^2_{+}}\setminus\Omega})\right)^4$.\\
\\
{\bf Lemma 4.3.} {\em For ${\mathfrak{F}}(x,y)$ defined by (\ref{Fxy}),
it holds that
\begin{equation}\label{eq-Fxy}
\frac{\partial {\mathfrak{F}}}{\partial x }(x,y)B+B\frac{\partial {\mathfrak{F}}}{\partial y
}(x,y)=0
\end{equation}
and
\begin{equation}\label{cond-Fxy}
{\mathfrak{F}}(x,0)={\mathcal{J}}(x),\ \
{\mathfrak{F}}(0,y)={\mathcal{L}}(y),
\end{equation}
where
\begin{equation}\label{L-L}
{\mathcal{J}}(x)=\overline{K^T\left(-\overline{P^T},0;\overline{Q^T}\right)(x,0)},\
\ {\mathcal{L}}(y)=K(P,0;Q)(y,0).
\end{equation}
Moreover, the following relation holds:
\begin{equation}\label{r-L-L}
{\mathcal{J}}(x)-B{\mathcal{J}}(x)B={\mathcal{L}}(x)-B{\mathcal{L}}(x)B.
\end{equation}}
 {\bf Proof.} For $y\leq x$, in view of (\ref{eq-K})$\sim$(\ref{xx-K}) in Lemma 2.1, we
see by integration by parts that
$$\ba{lll}&\ \ \ \displaystyle{\frac{\partial {\mathfrak{F}}}{\partial x }(x,y)B+B\frac{\partial {\mathfrak{F}}}{\partial y
}(x,y)}\\
\\&=\displaystyle{\Big\{ R(P,0)(y)\overline{
K^T_x\left(-\overline{P^T},0;\overline{Q^T}\right)(x,y)}}\\
\\
&\ \ \ \ \ \ \ \ \ \ \displaystyle{+\int_0^y K(P,0;Q)(y,s)\overline{
K^T_x\left(-\overline{P^T},0;\overline{Q^T}\right)(x,s)}{\rm
d}s\Big\}B }\\
\\&\ \ \ \ \ \ +B \Big\{ \displaystyle{R'(P,0)(y)\overline{ K^T\left(-\overline{P^T},0;\overline{Q^T}\right)(x,y)}}\\
\\&\ \ \ \ \ \ \ \ \ \ \ \ \ \displaystyle{+K(P,0;Q)(y,y)\overline{
K^T\left(-\overline{P^T},0;\overline{Q^T}\right)(x,y)} }\Big\}
\\
\\&\ \ \ \ \ \ +B \Big\{ \displaystyle{  R(P,0)(y)\overline{
K^T_y\left(-\overline{P^T},0;\overline{Q^T}\right)(x,y)}}\\
\\&\ \ \ \ \ \ \ \ \ \ \ \ \ \displaystyle{+\int_0^y K_y(P,0;Q)(y,s)\overline{
K^T\left(-\overline{P^T},0;\overline{Q^T}\right)(x,s)}{\rm d}s} \Big\}\\
\\
&=\left\{
BR'(P,0)(y)-R(P,0)(y)P(y)+BK(P,0;Q)(y,y)-K(P,0;Q)(y,y)B\right\}\\
\\
&\ \ \ \ \times \overline{
K^T\left(-\overline{P^T},0;\overline{Q^T}\right)(x,y)}+K(P,0;Q)(y,0)B\overline{
K^T\left(-\overline{P^T},0;\overline{Q^T}\right)(x,0)}
\\
\\&\ +\displaystyle{\int_0^y
\left\{BK_y(P,0;Q)(y,s)+K_s(P,0;Q)(y,s)B-K(P,0;Q)(y,s)P(s)\right\}}\\
\\& \ \ \ \ \
\times \overline{
K^T\left(-\overline{P^T},0;\overline{Q^T}\right)(x,s)}{\rm d}s \\
\\&=0,
\ea
$$
where we have made use of  the relation: $B=QB+BQ$. For the case $x\leq y$, the proof of (\ref{eq-Fxy}) is similar. On the other hand,
(\ref{cond-Fxy}) is obvious by (\ref{Fxy}).

Furthermore, it can be directly verified  that the unique solution of problem (\ref{eq-Fxy}) and
(\ref{cond-Fxy}) is
$${\mathfrak{F}}(x,y)=\left\{ \ba{lll} \displaystyle{\frac{1}{2}\{{\mathcal{J}}(x+y)+
{\mathcal{J}}(x-y)\}-\frac{1}{2}B\{{\mathcal{J}}(x+y)-
{\mathcal{J}}(x-y) \}B,\ y\leq x, }\\
\\
\displaystyle{\frac{1}{2}\{{\mathcal{L}}(x+y)+{\mathcal{L}}(y-x)\}-\frac{1}{2}B\{{\mathcal{L}}(x+y)
-{\mathcal{L}}(y-x)\}B,\ x\leq y.}\ea\right.
$$
Consequently, (\ref{r-L-L}) follows from the continuity of
${\mathfrak{F}}(x,y)$ at $x=y$.
 $\hfill\square$\\

Now we apply Lemma 4.3 to show\\
\\
{\bf Lemma 4.4.} {\em It holds that
$$\Theta_{{\mathcal{J}}}(\rho)Q=\Theta_{{\mathcal{J}}}(\rho)=
\widetilde{\Theta}_{\mathcal{L}}(\rho)=Q\widetilde{\Theta}_{\mathcal{L}}(\rho).
$$}
{\bf Proof.}  By (\ref{Theta}), we have
$$\Theta_{{\mathcal{J}}}(\rho)=\int_0^\infty {{\mathcal{J}}}(x) S(x,i\rho){\rm d}x,\
\widetilde{\Theta}_{\mathcal{L}}(\rho)=\int_0^\infty
\widetilde{S}(x,i\rho){\mathcal{L}}(x){\rm d}x,
$$
where $S(x,i\rho)$ and $\widetilde{S}(x,i\rho)$ are given by (\ref{SS}).
 Since $Q^2=Q$, it is sufficient to prove that for all $x\geq 0$
 $${{\mathcal{J}}}(x)Q=Q{\mathcal{L}}(x),\
 {{\mathcal{J}}}(x)BQ=-QB{\mathcal{L}}(x).
 $$
First, multiplying right (\ref{r-L-L}) by $Q$, we obtain by
$QB+BQ=B$ and ${\mathcal{L}}(x)Q={\mathcal{L}}(x)$ which follows from (\ref{K-p1p2-0}) and (\ref{L-L}) that
\begin{equation}\label{L=}
\{{{\mathcal{J}}}(x)-B{{\mathcal{J}}}(x)B\}Q=
{\mathcal{L}}(x)Q-B{\mathcal{L}}(x)(B-QB)={\mathcal{L}}(x).
\end{equation}
Second, since it follows from (\ref{K-p1p2*}) and (\ref{L-L}) that $Q{{\mathcal{J}}}(x)={{\mathcal{J}}}(x)$, we have
$QB{{\mathcal{J}}}(x)=(B-BQ){{\mathcal{J}}}(x)=0$.
Consequently, it follows from (\ref{L=}) that
$$Q{\mathcal{L}}(x)=Q\{{{\mathcal{J}}}(x)-B{{\mathcal{J}}}(x)B\}Q
={{\mathcal{J}}}(x)Q-QB{{\mathcal{J}}}(x)BQ={{\mathcal{J}}}(x)Q.
$$
On the other hand, multiplying left (\ref{L=}) by $B$, we have by
$B^2=E$ that
$$B{{\mathcal{J}}}(x)Q-{{\mathcal{J}}}(x)BQ=B{\mathcal{L}}(x)=(QB+BQ){\mathcal{L}}(x),
$$
that is,
$${{\mathcal{J}}}(x)BQ+QB{\mathcal{L}}(x)=B\{{{\mathcal{J}}}(x)Q-Q{\mathcal{L}}(x)\}=0.
$$
Thus the proof of Lemma 4.4 is completed.
$\hfill\square$\\
\\
{\bf Proof of Theorem 2.} Let
${\mathcal{L}}_\sigma(x)=\gamma_\sigma(x){\mathcal{L}}(x)$,
${{\mathcal{J}}}_\sigma(x)=\gamma_\sigma(x){{\mathcal{J}}}(x)$, where the scalar function $\gamma_\sigma(x)$ is defined by (\ref{gamma-sigma}). It is obvious that both
${\mathcal{L}}_\sigma$ and ${{\mathcal{J}}}_\sigma$ are continuously differentiable matrix-valued functions with compact support.
Then it follows easily from Lemma 4.4 that $\Theta_{{\mathcal{J}}_\sigma}(\rho)=\widetilde{\Theta}_{{\mathcal{L}}_\sigma}(\rho)$.
Hence, combining (\ref{eq-phi}), (\ref{eq-phi*}), (\ref{SS}), $Q{{\mathcal{J}}}(\cdot)={{\mathcal{J}}}(\cdot)$,
 ${\mathcal{L}}(\cdot)Q={\mathcal{L}}(\cdot)$ and Lemma 4.1, we conclude easily  that the  following matrix-valued function
$${\mathfrak{F}}_\sigma(x,y):=\frac{1}{\pi}\int_{-\infty}^\infty
S(y,i\rho)\Theta_{{\mathcal{J}}_\sigma}(\rho)\widetilde{S}(x,i\rho){\rm
d}\rho=\frac{1}{\pi}\int_{-\infty}^\infty
S(y,i\rho)\widetilde{\Theta}_{{\mathcal{L}}_\sigma}(\rho)\widetilde{S}(x,i\rho){\rm
d}\rho
$$
satisfies the equation
$$U_xB+BU_y=0,
$$
and the conditions
$$U(x,0)={{\mathcal{J}}}_\sigma(x),\
U(0,y)={\mathcal{L}}_\sigma(y)
$$
for all $x,y>0$. Therefore, if we define
${\mathfrak{F}}_\sigma(0,0)={\mathcal{L}}(0)={{\mathcal{J}}}(0)$,
then  ${\mathfrak{F}}_\sigma(x,y)={\mathfrak{F}}(x,y)$ in the
domain $0\leq x,y\leq\sigma$, since $\gamma_\sigma(x)\equiv 1$ on $[0,\sigma]$ and the two matrix-valued functions
satisfy the same boundary problem as that in Lemma 4.3. Moreover, if $f,g\in{\left({\mathbb{K}}^2_\sigma(0,\infty)\right)^4}$, then it follows from (\ref{fF}) and (\ref{gG}) that $F(x)=G(x)=0$ for
$x>\sigma$. Consequently, it follows from (\ref{Theta}), (\ref{fFgG}), Lemma 4.1 and Lemma 4.2 that
\begin{equation}\label{expan-1}
\ba{lll}&\ \ \ \displaystyle{\int_0^\infty f(x)g(x){\rm d}x}\\
\\
&=\displaystyle{\int_0^\infty F(x)G(x){\rm d}x+\int_0^\infty\int_0^\infty
F(y) {\mathfrak{F}}(x,y)G(x){\rm d}x{\rm d}y}\\
\\
&=\displaystyle{\int_0^\infty F(x)G(x){\rm d}x+\int_0^\infty\int_0^\infty
F(y) {\mathfrak{F}}_\sigma(x,y)G(x){\rm d}x{\rm d}y}\\
\\
&=\displaystyle{\frac{1}{\pi}\int_{-\infty}^\infty
\Theta_F(\rho)\{Q+\Theta_{{\mathcal{J}}_\sigma}(\rho)\}\widetilde{\Theta}_G(\rho){\rm
d}\rho}\\
\\
&=\displaystyle{\frac{1}{\pi}\int_{-\infty}^\infty
\Phi_f(\rho)\{Q+\widetilde{\Theta}_{\mathcal{L}_\sigma}(\rho)\}\widetilde{\Phi}_g(\rho){\rm
d}\rho}.\ea
\end{equation}
Now define
\begin{equation}\label{D}
D(\rho)=\lim_{\sigma\rightarrow \infty}\{Q+\Theta_{{\mathcal{J}}_\sigma}(\rho)\}=\lim_{\sigma\rightarrow \infty}\{Q+\widetilde{\Theta}_{\mathcal{L}_\sigma}(\rho)\}
\end{equation}
where the limits exist in the sense of convergence of distributions.  Indeed, by (\ref{Theta}) and (\ref{SS}) we see that both $\Theta_{{\mathcal{J}}_\sigma}(\rho)$ and $\widetilde{\Theta}_{\mathcal{L}_\sigma}(\rho)$
  are linear combination of the Fourier cosine and  sine transform of some matrix-valued function with compact support. Then it follows from the property of the Fourier transform (see e.g. Page 105 in \cite{mar}) that $\Theta_{{\mathcal{J}}_\sigma}(\rho)\rightarrow\Theta_{\mathcal{J}}(\rho)$ and $\widetilde{\Theta}_{\mathcal{L}_\sigma}(\rho)\rightarrow\widetilde{\Theta}_{\mathcal{L}}(\rho)$ as $ \sigma\rightarrow\infty$ in the sense of distributions, whence $D(\rho)\in{(Z')^4}$. Therefore, by the definition (\ref{Theta}) we see
$$D(\rho)=\displaystyle{\frac{1}{\pi}\{Q+\Theta_{{\mathcal{J}}}(\rho)\}=\frac{1}{\pi}
\{Q+\widetilde{\Theta}_{\mathcal{L}}(\rho)\}}.
$$
Thus we can prove the Marchenko-Parseval equality (\ref{parseval-2}) similarly to (\ref{parseval-1}). Moreover, if one lets $g(t)=\varsigma(t)E$ or $f(t)=\varsigma(t)E$ where
$\varsigma(t)$ is defined by (\ref{varsigma}), then he can prove
(\ref{expansion-2}) similarly to (\ref{expansion-1}).
$\hfill\square$\\
\\
\section*{Acknowledgements}

 The author would like to thank heartily Professor Vladimir Alexandrovich Marchenko and Professor
 Masahiro Yamamoto for their great comments and suggestion. This work was partially supported by the National Natural Science Foundation of China (Grant No. 11101390), the Fundamental Research Funds for the Central Universities and JSPS Fellowship P05297.


\begin{thebibliography}{00}
\bibitem{bck}
Bairamov E, Cakar O, Krall A M.  An eigenfunction expansion for a quadratic pencil of a Schr\"{o}dinger operator
with spectral singularitie. J Differential Equations, 1999, 151: 268-289

\bibitem{berez}
Berezanskii J M.  Expansions in Eigenfunctions of Selfadjoint Operators.
  Providence: AMS, 1968

\bibitem{frei01}
Freiling G, Yurko V. Inverse Sturm-Liouville Problems and Their
Applications. New York: Nova, 2001

\bibitem{frei12} Freiling G, Yurko V. Determination of singular differential pencils from the Weyl function. Adv Dyn Syst Appl,
 2012, 7: 171-193
\bibitem{gelf}
Gel'fand I M, Levitan B M. On the determination of a differential equation
from its spectral function. Amer Math Soc Transl (2), 1955, 1: 253-304


\bibitem{gb}
 Gesztesy F, Simon B. The xi function. Acta Math, 1996, 176: 49-71

\bibitem{gks}
Gohberg I, Kaashoek M A, Sakhnovich A L. Pseudo-canonical systems with rational Weyl functions:
explicit formulas and applications. J Differential Equations, 1998, 146: 375-398

\bibitem{gording}
G{\aa}rding L. On the asymptotic properties of the spectral
function belonging to a self-adjoint semi-bounded extension of an
elliptic differential operator. Kungl Fysiog S\"{a}llsk i Lund
F\"{o}rh, 1954, 24: 1-18


\bibitem{higdon}
 Higdon R L. Initial-boundary value problems for linear hyperbolic
  systems. SIAM Rev, 1986, 28: 177-217

\bibitem{hormander}
H\"{o}rmander L.  The spectral function of an elliptic operator.
 Acta Math, 1968, 121: 193-218

\bibitem{ike}
Ikebe T.  Eigenfunction expansions associated with the Schr\"{o}dinger operators and their
applications to scattering theory. Arch Rational Mech Anal, 1960, 5: 1-34

\bibitem{ks}
Killip R, Simon B.
Sum rules and spectral measures of Schr\"{o}dinger operators with $L^2$ potentials.
Ann  Math (2), 2009, 170: 739-782

\bibitem{kl}
Kirsten K, Loya P. Spectral functions for the Schr\"{o}dinger operator on $\Bbb R^+$ with a singular potential. J Math Phys, 2010, 51: 053512

\bibitem{kodaira}
Kodaira K. The eigenvalue problem for odinary differential
equations of the second order and Heisenberg's theory of
$S$-matrices. Amer J Math, 1949, 71: 921-945

\bibitem{kotani}
Kotani S, Matano H. Differential Equations and Eigenfunction Expansions (Japanese). Tokyo: Iwanami, 2006

\bibitem{lesch}
Lesch M, Malamud M M. The inverse spectral problem for first
order systems on the half line. In: Oper Theory Adv Appl
117. Basel: Birkh\"{a}user, 2000, 199-238

\bibitem{levinson}
 Levinson N.  The expansion theorem for singular self-adjoint linear
differential operators.  Ann Math (2), 1954, 59:
300-315

\bibitem{lev50}
Levitan B M.  Proof of the theorem on the expansion in
eigenfunctions of self-adjoint differential equations (Russian).
 Doklady Acad Nauk SSSR, 1950,73: 651-654

\bibitem{lev91}
Levitan B M, Sargsjan I S.  Sturm-Liouville and Dirac
Operators. Dordrecht: Kluwer, 1991

\bibitem{mar63}
Marchenko V A, Expansion in eigenfunctions of non-self-adjoint
singular differential operators of second order (Russian). Mat
Sb, 1960, 52: 739-788

\bibitem{mar}
Marchenko V A. Sturm-Liouville Operators and Applications.
 Basel: Birkh\"{a}user, 1986

\bibitem{me}
Mennicken R, M\"{o}ller M.  Non-Self-Adjoint Boundary
Eigenvalue Problems. Amsterdam: North-Holland, 2003

\bibitem{miyazaki}
Miyazaki Y.  Asymptotic behavior of spectral functions for elliptic
operators with non-smooth coefficients.  J Funct Anal, 2004, 214: 132-154

\bibitem{n06}
Ning W Q. An inverse spectral problem for a nonsymmetric differential operator: Reconstruction of eigenvalue
 problem.  J Math Anal Appl, 2007, 327: 1396-1419

\bibitem{n08}
Ning W Q. On stability of an inverse spectral problem for a
nonsymmetric differential operator.
  J Inverse Ill-Posed Probl, 2009, 17: 289-310

\bibitem{ny04}
Ning W Q, Yamamoto M. An inverse spectral problem for a nonsymmetric differential operator:
 Uniqueness and reconstruction formula.  Integr equat Oper Theory, 2006, 55: 273-304

\bibitem{ny08}
Ning W Q, Yamamoto M.  The Gel'fand-Levitan theory for
one-dimensional hyperbolic systems with impulsive inputs.
Inverse Problems, 2008, 24: 025004

\bibitem{sak}
Sakhnovich A.  Weyl functions, the inverse problem and special
solutions for the system auxiliary to the nonlinear
optics equation, Inverse Problems, 2008, 24: 025026

\bibitem{stone}
Stone M H.  Linear Transformations in Hilbert Space. Reprint of the 1932 original.
Providence: AMS, 1990


\bibitem{titchmarsh}
Titchmarsh E C.  Eigenfunction Expansions Associated
    with Second-Order Differential Equations. Oxford: Clarendon, 1962

\bibitem{tr}
Trooshin I.  On inverse scattering for nonsymmetric operators. Preprint

\bibitem{weyl}
Weyl H. \"{U}ber gew\"{o}hnliche Differentialgleichungen mit Singularit\"{a}ten und die zugeh\"{o}rigen
Entwicklungen willk\"{u}rlicher Funktionen (German).
 Math Ann, 1910, 68: 220-269


\bibitem{ya88}
Yamamoto M.  Inverse spectral problem for systems of ordinary
differential equations of first order, I. J Fac Sci Univ Tokyo Sect IA Math, 1988, 35: 519-546

\bibitem{yosida50}
Yosida K. On Titchmarsh-Kodaira formula concerning Weyl-Stone's
eigenfunction expansion.  Nagoya Math J, 1950, 1: 49-58

\bibitem{yurko92}
 Yurko V A.  Reconstruction of nonselfadjoint differential operators on the semi-axis from the Weyl matrix.
 Math USSR-Sb, 1992, 72: 413-438

\bibitem{yurko02}
Yurko V A. Method of Spectral Mappings in the Inverse Problem
Theory. Utrecht: VSP, 2002
\end{thebibliography}
\end{document}